\documentclass[11pt]{article}

\usepackage{amsmath,amssymb}
\usepackage{latexsym}
\usepackage{graphicx}
\usepackage{bm}

\newtheorem{thm}{Theorem}

\newtheorem{rem}{Remark}

\newenvironment{proof}{\begin{trivlist}
                       \item[]{\bf Proof.}
                       \hspace{0cm}}{\hfill $\Box$
                       \end{trivlist}}

\usepackage{hyperref}

\begin{document}
\title{An inverse problem for a heat equation with piecewise-constant 
thermal conductivity}

\author{N. S. Hoang$\dag$\footnotemark[1] \qquad A. G. Ramm$\dag$\footnotemark[2]
\\
\\
$\dag$Mathematics Department, Kansas State University,\\
Manhattan, KS 66506-2602, USA
}

\renewcommand{\thefootnote}{\fnsymbol{footnote}}
\footnotetext[1]{Email: nguyenhs@math.ksu.edu}
\footnotetext[2]{Corresponding author; Email: ramm@math.ksu.edu}

\date{}
\maketitle

\begin{abstract}
The governing equation is $u_t = (a(x)u_x)_x$, $0\le x\le 1$, $t>0$,
 $u(x,0)=0$, $u(0,t)=0$, $a(1)u'(1,t)=f(t)$.
 The extra data are $u(1,t)=g(t)$. It is assumed that
 $a(x)$ is a piecewise-constant function, and $f\not\equiv 0$. It is 
proved that
the function $a(x)$ is uniquely defined by the above data. No restrictions
on the number of discontinuity points of $a(x)$ and on their 
locations are made. The number of discontinuity points is finite,
but this number can be arbitrarily large.

If $a(x)\in C^2[0,1]$, then a uniqueness theorem has been established 
earlier
for multidimensional problem, $x\in \mathbb{R}^n, n>1$ (see MR1211417 
(94e:35004)) for the stationary problem with infinitely 
many boundary data.
The novel point in this work is the treatment of the discontinuous 
piecewise-constant function $a(x)$ and the proof of Property C for a pair
of the operators $\{\ell_1, \ell_2 \}$, where $\ell_j:= -\frac{d^2}{dx^2} 
+ k^2 q_j^2(x)$, $j=1,2$, and $q_j^2(x)>0$ are piecewise-constant
functions, and for 
the pair $\{L_1, L_2 \}$, where $L_ju:=-[a_j(x)u'(x)]'+\lambda u$,
$j=1,2$, and $a_j(x)>0$ are piecewise-constant functions.
Property C stands for completeness of the set of products of solutions
of homogeneous differential equations (see MR1759536 (2001f:34048))  

{\bf Keywords}: inverse problems, heat equation, Property C, piecewise-
constant thermal conductivity.

{\bf MSC}:  35R30, 74J25, 34E05.
\end{abstract}

\section{Introduction}

Let 
\begin{align}
\label{eq1}
\dot{u} = (a(x)u')', \quad 0\le x\le 1,\quad t>0,\quad 
u':=\frac{\partial u}{\partial x},\quad
\dot{u}:=\frac{\partial u}{\partial t},\\
\label{eq2}
u(x,0) = 0,\quad u(0,t) = 0,\quad a(1)u'(1,t) = f(t)\not\equiv 0,\\
\label{eq3}
u(1,t) = g(t).
\end{align}
Problem \eqref{eq1}--\eqref{eq2} describes the heat transfer in a rod, $a(x)$
is the heat conductivity, $a(1)u'(1,t)$ is the heat flux, $g(t)$ 
is the measurement, the extra data. 

The inverse problem (IP) is:

\textbf{IP:} {\it Given $f(t)$ and $g(t)$ for all $t>0$, find $a(x)$}.

{\bf Assumption A:} {\it $a(x)$ is a piecewise-constant function, 
$a(x)=a_j$, 
$x_j\le x\le x_{j+1}$,
$x_1=0$, $x_{n+1}=1$, $0<c_0\le a_j\le c_1$, $1\le j\le n$.}

This assumption holds throughout the paper and is not repeated. The set 
of piecewise-constant functions with finitely many discontinuity points is
denoted by $\Pi $.

If $a(x)\in C^2$, then the uniqueness of the solution to 
some multidimensional inverse problems has been proved in \cite{R278}
(see also \cite{R470}).
Problem \eqref{eq1}--\eqref{eq3} with $a(x)\in C^2([0,1])$ has been studied 
in \cite{R422}, \cite{R438}. The treatment of discontinuous piecewise-constant
$a(x)$ is of interest in applications.

In \cite{GH} equation \eqref{eq1} with the conditions $u(0,t)=u(1,t)=0$,
$u(x,0)=g(x)$, was studied, and the measured (extra) data were the values
$u(\xi_m,t)$, $\forall t>0$, $1\le m\le M$,
$0\le \xi_m\le 1$, where $M=3n$, and $n$ is the number of the discontinuity points
of $a(x)$. It was assumed in \cite{GH} that $\min_j|x_j-x_{j+1}|$
is not too small. Under these assumptions the uniqueness theorem for the IP
was proved in \cite{GH}, and an algorithm for finding $a(x)$ was proposed.
The stability of this algorithm with respect to perturbations of the data was not studied in \cite{GH}.

In our paper the extra data \eqref{eq3} consists of measurement, taken 
at one point,  rather than at $3n$ points, and we impose
no restrictions on $\min_j|x_j - x_{j+1}|$. Under these assumptions, 
which are much weaker than in \cite{GH}, we prove the 
uniqueness of the solution to IP. 

One of our main results is

\begin{thm}
\label{theorem1}
The IP has at most one solution.
\end{thm}

\begin{rem}
The IP is ill-posed: small variations of the data $\{f(t), g(t)\}$
in the $C(0,\infty)-$norm may lead to large variations of the coefficient 
$a(x)$, or may lead to a problem which has no solutions. We assumed that 
the data are known for all $t>0$. If one assumes that $f(t)=0$ for $t>T$,
where $T>0$ is an arbitrary fixed number, then the solution
$u(x,t)$ is an analytic function of $t$ in the region $t>T$.
Therefore the  data $\{f(t), g(t)\}$, known in the interval 
$[0,T+\epsilon)$, where $\epsilon>0$ is an arbitrary small fixed number, 
determine uniquely the data for all $t>0$. Thus, if 
$f(t)=0$ for $t>T$, then the uniqueness theorem for the solution to IP 
remains valid if the data are known for $t\in [0,T+\epsilon)$.
\end{rem}

Let us formulate IP in an equivalent form.

Take the Laplace transform of the equation \eqref{eq1}--\eqref{eq3},
denote 
$$v(x,\lambda):=Lu:=\int_0^\infty e^{-\lambda t}u(x,t)dt,$$
 and get:
\begin{align}
\label{eq4}
\lambda v - (a(x)v')'= 0,\quad 0\le x\le 1,\quad v(0,\lambda) = 0,\\
\label{eq5}
a(1)v'(1,\lambda) = F(\lambda),\qquad v(1,\lambda) = G(\lambda),
\end{align}
where $F:=Lf$ and $G:=Lg$.

The IP can be reformulated as follows:

\textbf{IP}: {\it Given $F(\lambda)$ and $G(\lambda)$ for all 
$\lambda>0$, find $a(x)$. }

Let us transform equation \eqref{eq4}-\eqref{eq5} to yet another equivalent form.

Let $a(x)v':=\psi$. Then \eqref{eq4}--\eqref{eq5} can be replaced by the following problem
\begin{align}
\label{eq6}
-\psi'' +\lambda a^{-1}(x)\psi = 0,\quad \psi(1,\lambda)=F(\lambda),\quad 
\psi'(0,\lambda) = 0,\\
\label{eq7}
\psi'(1,\lambda) = \lambda G(\lambda).
\end{align}

The IP can be reformulated as follows:

\textbf{IP} Given $G(\lambda)$ and $F(\lambda)$, find $a^{-1}(x):=q^2(x)$.

Let
\begin{equation}
\label{eq8}
\ell \psi := -\psi'' + k^2 q^2(x)\psi  = 0,\quad \lambda:=k^2,\quad 
q^2(x):=a^{-1}(x),\quad c_1^{-1}\le q^2(x) \le c_0^{-1}. 
\end{equation}
Consider the following problems:
\begin{equation}
\label{eq9}
\ell_j\psi_j = 0,\quad \ell_j:= -\frac{d^2}{dx^2} + k^2 q_j^2(x),\quad \psi_j'(0,k)=0,
\quad \psi_j(0,k)=1,\quad j=1,2.
\end{equation}

Our second main result is

\begin{thm}
\label{theorem2}
The sets $\{\psi_1(x,k)\psi_2(x,k)\}_{\forall k\ge 0}$ 
and $\{v'_1(x,\lambda)v'_2(x,\lambda)\}_{\forall \lambda\ge 0}$, 
$k:=\lambda^{1/2},$ are dense 
in the set $\Pi$
of piecewise-constant functions on $[0,1]$.
\end{thm}

\begin{rem}
Theorem~\ref{theorem2} says that if $h(x)\in\Pi$ and 
\begin{equation}
\label{eq10}
\int_0^1 h(x)\psi_1(x,k)\psi_2(x,k)dx = 0,\qquad \forall k>0,
\end{equation} 
then $h=0$. Similar conclusion holds if $\psi_j(x,k)$ is replaced by 
$v'_j(x,\lambda)$ in \eqref{eq10}.
Such a property of the pair of the operators $\{\ell_1,\ell_2\}$ is called
Property C (\cite{R470}, \cite{R402}). 

Clearly if the set $\{\psi_1(x,k)\psi_2(x,k)\}_{\forall k\ge 0}$
is dense in the set $\Pi$, then the set of products 
$\{v'_1(x,\lambda)v'_2(x,\lambda)\}_{\forall \lambda\ge 0}$
is dense in the set $\Pi$.
\end{rem}

In Section 2 proofs are given.

\section{Proofs}

\subsection{Proof of Theorem~\ref{theorem1}}
\begin{proof}
We prove this Theorem for the problem \eqref{eq4}--\eqref{eq5}.
Suppose there are $v_j$ and $a_j\in \Pi$, $j=1,2$, which solve
problem \eqref{eq4}--\eqref{eq5}, and let $w:=v_1-v_2$.
Then
\begin{align}
\label{eq11}
\lambda w - (a_1w')' = (pv_2')',\qquad p:=a_1(x)-a_2(x),\\
\label{eq12}
w(0,\lambda) = 0,\quad w(1,\lambda)=0, \quad 
a_1v_1'(1,\lambda)=a_2v_2'(1,\lambda).
\end{align}
Multiply \eqref{eq11} by $v_1$, a solution to equation \eqref{eq4} with
$a=a_1$, and
integrate over $[0,1]$, and then by parts, to get
\begin{equation}
\label{eq13}
\int_0^1 p(x)v_2'v_1'dx = pv_2'v_1\bigg{|}_0^1
+ a_1w'v_1\bigg{|}_0^1 - a_1w v_1'\bigg{|}_0^1 = 0,\qquad \forall 
\lambda>0,\quad \lambda=k^2,\,\, k>0,
\end{equation}
where we have used the conditions $w(0,\lambda)=w(1,\lambda)=0$
and $a_1(1)v_1'(1,\lambda)=a_2(1)v_2'(1,\lambda)$.
Note that $v_2(x,\lambda)$ can be considered as an arbitrary 
solution to equation  \eqref{eq4},  up to
a constant factor.
The set $\{v'_1(x,\lambda)v'_2(x,\lambda)\}$ is dense in $\Pi$ by 
Theorem~\ref{theorem2}.
Since $a_1(x)-a_2(x):=p(x)\in \Pi$,  it follows 
from \eqref{eq13} that
$p(x)=0$. So $a_1=a_2$. Theorem~\ref{theorem1} is proved.
\end{proof}

\subsection{Proof of Theorem~\ref{theorem2}}

\begin{proof}
Let us prove completeness of the set of products 
$\{\psi_1(x,k)\psi_2(x,k)\}_{\forall k\ge 0}$.
Assume that $h\in \Pi$ and \eqref{eq10} holds. The function $\psi_j(x,k)$, 
$j=1,2$,
are entire functions of $k$. This follows from the integral 
equation for $\psi_j$,
which is an immediate consequence of equations \eqref{eq8}--\eqref{eq9}:
\begin{equation}
\label{eq14}
\psi_j(x,k) = 1 + k^2\int_0^x(x-s)q_j^2(s)\psi_j(s,k)ds,\qquad x\ge 0,
\quad j=1,2.
\end{equation}
Equation \eqref{eq14} implies that for any fixed $k$ one has 
$\psi_i(x):=\psi_j(x,k)\geq 1$,    
$\forall x\in [0,1]$, $j=1,2$, that $\psi'_j(x,k)\geq 0$,
 $\psi''_j(x,k)\geq 0$, and $\frac {\partial^m \psi_j(x,k)}
{\partial k^m}\geq 0$ for all $m=0,1,2,.....$. Consequently,
$\psi_j(x)$, $j=1,2$, are convex functions of $x$ on the semiaxis $x>0$.
Since $\psi_j(x,k)$, $j=1,2$, are positive, it follows from \eqref{eq14} 
that $\psi_j(x,k)$, $j=1,2$, are 
increasing functions with respect to both $x$ and $k$. So we have
\begin{equation}
\label{eq15}
\psi_j(x,k)>0,\quad \psi_j'(x,k)>0,\quad \psi_j''(x,k)>0,\qquad \forall k>0,\quad j=1,2.
\end{equation}

Assume $0< x_{11}< x_{12}<\cdots < x_{1N_1}<1$ and $0< x_{21}< 
x_{22}<\cdots < x_{2N_2}<1$
are discontinuity points of $a_1(x)$ and $a_2(x)$, respectively.

{\it To derive from \eqref{eq10} that $h=0$ it is sufficient to prove 
that} $h(x)=0,\forall x\in [x_0,1]$, where $x_0:=\max(x_{1N_1}, 
x_{2N_2})$,
because then one can prove similarly, in finitely many steps,
that $h=0$ on the whole interval $[0,1]$ using the assumption $h\in \Pi$.
We have
\begin{equation}
\label{eq16}
\psi''_j(x,k) = k^2q_{jN_j}^2(x)\psi_j(x,k),\quad \forall k>0,\quad 
\forall 
x\in 
[x_{0},1],
\end{equation}
where $q_{jN_j}$ is the value of $q_j$ on the interval $[x_0,1]$.
From \eqref{eq16} one gets
\begin{equation}
\label{eq17}
\psi_j(x,k) = a_j(k)e^{kq_{jN_j}(x-x_0)} + 
b_j(k)e^{-kq_{jN_j}(x-x_0)},\qquad \forall k\ge 0,\quad j=1,2.
\end{equation}
It follows from \eqref{eq15} and \eqref{eq17} that
\begin{equation}
\label{eq18}
\psi_j(x_0,k) = a_j(k) + b_j(k)>0,\quad 
\psi_j'(x_0,k) = kq_{jN_j}[a_j(k) - b_j(k)] \ge 0,
\end{equation}
and 
\begin{equation}
\label{eq20}
2a_j(k)=\psi_j(x_0,k) + \frac 
{\psi_j'(x_0,k)}{kq_{jN_j}}>\psi_j(x_0,k).
\end{equation}
This implies
\begin{equation}
\label{eq19}
a_j(k)\ge |b_j(k)|\geq 0,\qquad \forall k>0,\quad j=1,2.
\end{equation} 

Since $h\in \Pi$, one may assume without loss of generality that 
\begin{equation}
\label{eq24}
h(x)=C\ge 0,\qquad \forall x\in [x_0,1].
\end{equation}
It follows from \eqref{eq10} that
\begin{equation}
\label{eq25}
-\int_0^{x_0}\psi_1(x,k)\psi_2(x,k)h(x)dx = \int_{x_0}^1\psi_1(x,k)\psi_2(x,k)h(x)dx,\qquad \forall k>0.
\end{equation}
From \eqref{eq15}, \eqref{eq17} and \eqref{eq19}, one gets
\begin{equation}
\label{eq26}
1\le \psi_j(x,k)\le \psi_j(x_0,k)< 2a_j(k),\quad 0\leq x\leq x_0,\qquad 
\forall k>0, \quad j=1,2.
\end{equation}
Therefore,
\begin{equation}
\label{eq27}
\bigg{|} \int_0^{x_0}\psi_1(x,k)\psi_2(x,k)h(x)dx\bigg{|} 
\le 4a_1(k)a_2(k)\int_0^{x_0}|h(x)|dx.
\end{equation}
From \eqref{eq19}, \eqref{eq17} and \eqref{eq15} one obtains
\begin{equation}
\label{eq28}
\psi_j(x,k)\geq a_j(k)[e^{kq_{jN_j}(x-x_0)} - e^{-kq_{jN_j}(x-x_0)}],\qquad 
x\in 
[x_0,1],\quad j=1,2.
\end{equation}
Take an arbitrary $y\in (x_0,1)$ and fix it. One has
$\psi_j(x,k) \ge \psi_j(y,k),\,
\forall x\in [y,1]$. Therefore,
\begin{equation}
\label{eq29}
\int_{x_0}^1 \psi_1(x,k)\psi_2(x,k)h(x)dx\geq
C(1-y)\psi_1(y,k)\psi_2(y,k),\qquad \forall k>0.
\end{equation}
This, \eqref{eq26}, \eqref{eq25}, and \eqref{eq27} imply the following inequalities:
\begin{equation}
\label{eq30}
\infty> 4\int_0^{x_0}|h(x)|dx \ge
C(1-y)\frac{\psi_1(y,k)\psi_2(y,k)}{a_1(k)a_2(k)},\qquad \forall
k>0. \end{equation}
It follows from \eqref{eq28} that
\begin{equation}
\label{eq31}
\lim_{k\to\infty} \frac{\psi_j(y,k)}{a_j(k)} = \infty.
\end{equation}
Let $k\to \infty$ in \eqref{eq30} and use \eqref{eq31} to conclude
that $C=0$ and, therefore, $h(x)=0$ for $x\in [x_0,1]$.
Similarly one proves that $h(x)=0$ for all $x\in [0,1]$.

Theorem~\ref{theorem2} is proved.
\end{proof}

\end{document}